\documentclass{article}

\usepackage{graphicx}
\usepackage{amsmath}
\usepackage{amsfonts}
\usepackage{amssymb}
\usepackage[french]{babel}

\newtheorem{theorem}{Th\'eor\`eme}
\newtheorem{lemma}{Lemme}
\newtheorem{proposition}{Proposition}

\newtheorem{definition}{D\'efinition\rm}

\newtheorem{remafin}{Remarque finale}

\newenvironment{proof}[1][D\'emonstration]{\textbf{#1.} }{\ \rule{0.5em}{0.5em}}

\begin{document}

\title{Les alg\`ebres de Lie r\'esolubles rigides r\'eelles ne sont pas n\'ecessairement compl\`etement r\'esolubles
\footnote{Les auteurs ont \'et\'e support\'es par le projet de recherche PR1/05-13283 de la UCM.}}

\author{J. M. Ancochea Berm\'udez\footnote{e-mail: ancochea@mat.ucm.es}, R. Campoamor-Stursberg\footnote{e-mail: rutwig@mat.ucm.es}\\ and L. Garc\'{\i}a Vergnolle\footnote{e-mail: lucigarcia@mat.ucm.es}\\
 \\
Dpto. Geometr\'{\i}a y Topolog\'{\i}a, Facultad CC. Matem\'aticas U.C.M.\\Plaza de Ciencias 3, E-28040 Madrid\\ }

\date{}

\maketitle

\begin{abstract}
On montre qu'une alg\`ebre de Lie r\'esoluble rigide r\'eelle n'est pas n\'ecessairement compl\`etement r\'esoluble. On  construit un exemple $\frak{n}\oplus\frak{t}$ de dimension minimale dont le tore ext\'erieur $\frak{t}$ n'est pas form\'e par des d\'erivations $ad$-semi-simples sur $\mathbb{R}$. Nous \'etudions les formes r\'eelles des nilradicaux des alg\'ebres de r\'esolubles rigides en dimension $n\leq 7$ et donnons la classification des alg\`ebres r\'esolubles rigides sur $\mathbb{R}$  en dimension 8.

\end{abstract}

\section{Preliminaries}

Rappellons qu'une alg\`ebre de Lie $\frak{g}$ sur un corps $\mathbb{K}$ de caract\'eristique 0 est dite rigide lorsque son orbite dans le sch\'ema $L_{n}$ donn\'e par les r\'elations de Jacobi sur la famille des lois d'alg\`ebres de Lie sur $\mathbb{K}^{n}$ est ouverte sous l'action du groupe $GL(n,\mathbb{K})$. Si $H^{2}(\frak{g},\frak{g})$ est le deuxi\`eme groupe de cohomologie \`a valeurs dans le module adjoint, le crit\`ere de rigidit\`e de Nijenhuis et Richardson \cite{NR} donne une condition suffisante: si $H^{2}(\frak{g},\frak{g})=0$, alors l'alg\`ebre $\frak{g}$ est rigide dans $L_{n}$. Dans \cite{Ca}, on montre que sur un corps $\mathbb{K}$ alg\'ebriquement clos de caract\'eristique 0, toute alg\`ebre de Lie r\'esoluble rigide v\'erifie la d\'ecomposition suivante:
\begin{equation*}
\frak{g}=\frak{n}\oplus\frak{t},
\end{equation*}
o\`u $\frak{n}$ est le nilradical de $\frak{g}$ et $\frak{t}$ est une sous-alg\`ebre ab\'elienne dont les \'el\'ements sont $ad_{\frak{g}}$-semi-simples, appel\'e tore ext\'erieur de $\frak{g}$. En utilisant cette propri\'et\'e, on peut d\'evelopper une th\'eorie de racines qui permet classifier les alg\'ebres de Lie r\'esolubles rigides complexes par rapport au spectre du tore ext\'erieur.

\begin{definition}
Un \'el\'ement $X\in\frak{t}$ est dit r\'egulier si la dimension du noyau ${\rm ker }\, ad(X)$ est minimale.
\end{definition}

Comme $ad(X)$ est diagonalisable et le nilradical est invariant par cet op\'erateur, on peut trouver une base $\left\{X_{1},..,X_{n}=X\right\}$ de vecteurs propres telle que $\left\{X_{1},..,X_{p+q}\right\}$ soit une base de $\frak{n}$, $\left\{X_{p+q+1},..,X_{n}=X\right\}$ soit une base de $\frak{t}$ et $\left\{X_{p+1},..,X_{n}=X\right\}$ une base de ${\rm ker}\;ad(X)$. Alors on consid\`ere le syst\`eme lin\'eaire de racines $S$ \`a $\dim\frak{g}-1$ variables $x_{i}$ dont les \'equations sont  $x_{i}+x_{j}=x_{k}$ si la composante du vecteur $[X_{i},X_{j}]$ sur $X_{k}$ est non nulle. De cette fa\c{c}on, on obtient un crit\`ere de rigidit\'e:

\begin{proposition}
{\rm \cite{A2}} Si $\frak{g}$ est rigide, alors pour tout vecteur r\'egulier $X$ on a ${\rm rank}(S)=\dim\frak{n}-1$.
\end{proposition}

En utilisant ce syst\`eme, on a classifi\'e les alg\`ebres de Lie r\'esolubles rigides complexes jusqu'\`a dimension huit \cite{A2}. De plus, cette propri\'et\'e lin\'eaire implique une autre d'ordre structurale:

\begin{definition}
Une alg\`ebre de Lie $\frak{r}$ sur le corps $\mathbb{K}$ est dite compl\`etement r\'esoluble si pout tout $X\in\frak{r}$ les valeurs propres de $ad X$ sont des \'el\'ements de $ \mathbb{K}$.
\end{definition}

En particulier, cette structure implique l'existance d'une suite
\begin{equation*}
\frak{r}=I_{n}\supset I_{n-1}\supset...\supset I_{0}=0,
\end{equation*}
o\`u les $I_{k}$ sont des id\'eaux de codimension 1 dans $I_{k+1}$.

\begin{proposition}
Soit $\frak{r}$ une alg\`ebre de Lie r\'esoluble rigide complexe.
Alors $\frak{r}$ est compl\`etement r\'esoluble.
\end{proposition}

\begin{proof}
D'apr\`es la proposition 1, pour toute alg\`ebre de Lie r\'esoluble complexe rigide il existe une base $\left\{X_{1},..,X_{r}\right\}$ du tore ext\'erieur $\frak{t}$ telle que les valeurs propres des op\'erateurs adjointes sont
des enti\`eres. En effet, \'etant donn\'ee une base $\left\{Y_{1},..,Y_{r}\right\}$ de $\frak{t}$, les valeurs propres des op\'erateurs $ad(Y_{j})$ sont des solutions du syst\`eme $S(X)$. Comme celui-ci est \`a coefficient entiers, il existe une base de solutions rationelles, et pourtant on peut trouver une base de solutions enti\`eres.
\end{proof}

\medskip

Dans \cite{A3} on a montr\'e l'existence des alg\`ebres de Lie rigides complexes non-r\'eelles. D'apr\`es le th\'eor\`eme du rang, \c{c}a implique que le tenseur de structure du nilradical n'est pas r\'eel, c'est-\`a-dire, le nilradical n'admet pas des formes r\'eelles\footnote{Rappellons qu'une forme r\'eelle d'une alg\`ebre de Lie complexe $\frak{g}$ est une alg\`ebre r\'eelle $\frak{g}^{\prime}$ telle que $\frak{g}^{\prime}\otimes\mathbb{C}\simeq \frak{g}$.}. On appellera une alg\`ebre de Lie complexe sans formes r\'eelles purement complexe.

\begin{proposition}
Toute alg\`ebre de Lie r\'esoluble rigide complexe $\frak{g}=\frak{n}\oplus\frak{t}$ dont le nilradical n'est pas purement complexe admet des formes r\'eelles.
\end{proposition}

La preuve est imm\'ediate. En effet, si le nilradical n'est pas purement complexe, alors les constantes de structure de $\frak{n}$ sont r\'eelles. Par le th\'eor\`eme du rang, on peut trouver des g\'en\'erateurs du tore ext\'erieur ayant des valeurs propres enti\`eres. Alors la restriction par scalaires donne une alg\`ebre de Lie r\'eelle, appel\'ee la forme r\'eelle normale. De plus, toute autre forme r\'eelle $(\frak{n}^{\prime}\oplus\frak{k})$ poss\`ede au moins un g\'en\'erateur $X\in\frak{k}$ tel que $ad(X)$ n'est pas diagonalisable sur $\mathbb{R}$, mais $X\in\frak{t}$ vu comme g\'en\'erateur de l'alg\`ebre complexifi\'ee.

\section{L'exemple}

Dans ce paragraphe on montre, par construction d'un exemple, que pour les alg\`ebres de Lie r\'eelles r\'esolubles rigides la propri\'et\'e de r\'esolubilit\'e compl\`ete ne reste pas valable si la forme r\'eelle n'est pas la normale. L'exemple construit est aussi de dimension minimale.

\smallskip

Consid\'{e}rons les lois d'alg\`{e}bres de Lie nilpotentes $\frak{g}_{1}=(\mathbb{R}^{6},\mu _{1})$ et $\frak{g}_{2}=(\mathbb{R}^{6},\mu _{2})$ d\'{e}finies dans une base $\left\{ X_{1},..,X_{6}\right\} $ par les crochets:

\begin{eqnarray*}
&&
\begin{tabular}{llll}
$\mu _{1}\left( X_{1},X_{i}\right) =X_{i+1}$ & $\left( 2\leq i\leq 4\right) ,
$ & $\mu _{1}\left( X_{3},X_{2}\right) =X_{6},$ & $\mu _{1}\left(
X_{6},X_{2}\right) =X_{5},$%
\end{tabular}
\\
&&
\begin{tabular}{llll}
$\mu _{2}\left( X_{1},X_{i}\right) =X_{i+1}$ & $\left( 2\leq i\leq 4\right) ,
$ & $\mu _{2}\left( X_{3},X_{2}\right) =-X_{6},$ & $\mu _{2}\left(
X_{6},X_{2}\right) =X_{5},$%
\end{tabular}
\end{eqnarray*}
Les alg\`{e}bres ne sont pas isomorphes, et
correspondent aux classes d'isomorphisme $N_{6,6}$ et $N_{6,7}$ de la liste \cite{Ce}. Si on consid\`ere le produit tensoriel, on a $\frak{g}_{1}\otimes\mathbb{C}\simeq \frak{g}_{2}\otimes\mathbb{C}$.

\begin{lemma}
L'alg\`{e}bre $\frak{g}_{1}=\left( \mathbb{R}^{6},\mu _{1}\right) $ admet une d\'{e}%
rivation diagonalisable $f_{1}^{1}$ et une d\'{e}rivation non-nilpotente $%
f_{1}^{2}$ donn\'{e}es respectivement par
\begin{eqnarray*}
&&
\begin{tabular}{lll}
$f_{1}^{1}\left( X_{i}\right) =X_{i}\;\left( i=1,2\right) ,$ & $%
f_{1}^{1}\left( X_{i}\right) =(i-1)X_{i}\;\left( i=3,4,5\right) ,$ & $%
f_{1}^{1}\left( X_{6}\right) =3X_{6}.$%
\end{tabular}
\\
&&
\begin{tabular}{llll}
$f_{1}^{2}\left( X_{1}\right) =X_{2},$ & $f_{1}^{2}\left( X_{2}\right)
=-X_{1}$ & $f_{1}^{2}\left( X_{4}\right) =-X_{6}.$ & $f_{1}^{2}\left(
X_{6}\right) =X_{4}.$%
\end{tabular}
\end{eqnarray*}
sur la base $\left\{ X_{1},..,X_{6}\right\}$. De plus, $f_{1}^{1}\circ
f_{1}^{2}=f_{1}^{2}\circ f_{1}^{1}$.
\end{lemma}

C'est facile de voir que toute d\'{e}rivation non nilpotente de l'alg%
\`{e}bre $\left( \mathbb{R}^{6},\mu _{1}\right) $ fait intervenir une combinaison lin\'eaire des d\'erivations $f_{1}^{1}$ et $f_{1}^{2}$. Comme le polyn\^{o}me caract\'{e}ristique de $f_{1}^{2}$ sur la
base donn\'{e}e est $P\left( T\right) =T^{2}\left( T^{2}+1\right) ^{2}$,
cette d\'{e}rivation n'est pas diagonalisable sur $\mathbb{R}$. De cette fa\c{c}on, on obtient une alg\`ebre abelienne $\frak{t}\in Der(\frak{g}_{1})$ constitu\'ee par d\'erivations non nilpotentes, mais dont seulement une est diagonalisable.

\begin{lemma}
L'alg\`{e}bre $\frak{g}_{2}=\left( \mathbb{R}^{6},\mu _{2}\right) $ admet un tore
exterieur $\frak{t}$ de rang 2.
\end{lemma}

Toute d\'{e}rivation non nilpotente de l'alg%
\`{e}bre $\left( \mathbb{R}^{6},\mu _{2}\right) $ fait intervenir une combinaison lin\'eaire des d\'erivations
\begin{eqnarray*}
&&
\begin{tabular}{lll}
$g_{1}^{1}\left( X_{i}\right) =X_{i}\;\left( i=1,2\right) ,$ & $%
g_{1}^{1}\left( X_{i}\right) =(i-1)X_{i}\;\left( i=3,4,5\right) ,$ & $%
g_{1}^{1}\left( X_{6}\right) =3X_{6}.$%
\end{tabular}
\\
&&
\begin{tabular}{llll}
$g_{1}^{2}\left( X_{1}\right) =X_{2},$ & $g_{1}^{2}\left( X_{2}\right) =X_{1}
$ & $g_{1}^{2}\left( X_{4}\right) =X_{6}.$ & $g_{1}^{2}\left( X_{6}\right)
=X_{4}.$%
\end{tabular}
\end{eqnarray*}
C'est imm\'{e}diat que les d\'{e}rivations sont diagonalisables, et sur la
base de vecteurs propres
\begin{eqnarray*}
\left\{ {X_{1}^{\prime }=\frac{1}{2}\left( X_{1}+X_{2}\right)
,X_{2}^{\prime
}=\frac{1}{2}(X_{2}-X_{1}),X_{3}^{\prime }=\frac{1}{2}X_{3},X_{4}^{\prime }=%
\frac{1}{4}(X_{4}-X_{6}),}\right. \\
\left. {X_{5}^{\prime }=\frac{1}{4}(X_{6}+X_{4}),X_{6}^{%
\prime }=\frac{1}{4}X_{5}}\right\}
\end{eqnarray*}
on obtient les crochets
\begin{equation*}
\begin{tabular}{llll}
$\mu _{2}\left( X_{1}^{\prime },X_{i}^{\prime }\right) =X_{i+1}^{\prime }$ &
$\left( i=2,3,5\right) ,$ & $\mu _{2}\left( X_{3}^{\prime },X_{2}^{\prime
}\right) =-X_{5}^{\prime },$ & $\mu _{2}\left( X_{2}^{\prime },X_{4}^{\prime
}\right) =X_{6}^{\prime }.$%
\end{tabular}
\end{equation*}

Les valeurs propres des d\'{e}rivations $g_{1}^{1}$ et $g_{1}^{2}$ sont,
respectivement
\begin{eqnarray*}
\sigma \left( g_{1}^{1}\right)  &=&\left( 1,1,2,3,3,4\right)  \\
\sigma \left( g_{1}^{2}\right)  &=&\left( -1,1,0,-1,1,0\right) .
\end{eqnarray*}
Pour commodit\'{e}, on consid\`{e}re les d\'{e}rivations $f_{2}^{1}=\frac{1}{2%
}\left( 3g_{1}^{1}+g_{1}^{2}\right) $ et $f_{2}^{2}=\frac{1}{2}\left(
g_{1}^{1}+g_{1}^{2}\right) $. Pour cette base du tore $\frak{t}$, les
valeurs propres sont:
\begin{eqnarray*}
&&
\begin{tabular}{lll}
$f_{2}^{1}\left( X_{i}^{\prime }\right) =iX_{i}^{\prime }$ & $\left( 1\leq
i\leq 6\right) ,$ &
\end{tabular}
\\
&&
\begin{tabular}{llll}
$f_{2}^{2}\left( X_{i}^{\prime }\right) =X_{i}^{\prime }$ & $\left( 2\leq
i\leq 4\right) ,$ & $f_{2}^{2}\left( X_{i}^{\prime }\right) =2X_{i}^{\prime }
$ & $\left( i=5,6\right) .$%
\end{tabular}
\end{eqnarray*}

\bigskip

Soient $\widehat{\frak{g}}_{1}=\left( \mathbb{R}^{8},\widehat{\mu }_{1}\right) $ et $\widehat{\frak{g}}_{2}=\left( \mathbb{%
R}^{8},\widehat{\mu }_{2}\right) $ les alg\`{e}bres r\'{e}solubles non
nilpotentes d\'{e}finies par
\begin{equation*}
\begin{tabular}{llll}
$\widehat{\mu }_{1}\left( X_{1},X_{i}\right) =X_{i+1}$ & $\left( 2\leq i\leq
4\right) ,$ & $\widehat{\mu }_{1}\left( X_{3},X_{2}\right) =X_{6},$ & $%
\widehat{\mu }_{1}\left( X_{6},X_{2}\right) =X_{5},$ \\
$\widehat{\mu }_{1}\left( X_{7},X_{i}\right) =X_{i}$ & $\left( i=1,2\right) ,
$ & $\widehat{\mu }_{1}\left( X_{7},X_{3}\right) =2X_{3},$ & $\widehat{\mu }%
_{1}\left( X_{7},X_{5}\right) =4X_{5},$ \\
$\widehat{\mu }_{1}\left( X_{7},X_{i}\right) =3X_{i}$ & $\left( i=4,6\right)
,$ & $\widehat{\mu }_{1}\left( X_{8},X_{1}\right) =X_{2},$ & $\widehat{\mu }%
_{1}\left( X_{8},X_{2}\right) =-X_{1},$ \\
$\widehat{\mu }_{1}\left( X_{8},X_{4}\right) =-X_{6},$ & $\widehat{\mu }%
_{1}\left( X_{8},X_{6}\right) =X_{4}.$ &  &
\end{tabular}
\end{equation*}
\begin{equation*}
\begin{tabular}{llll}
$\widehat{\mu }_{2}\left( X_{1},X_{i}\right) =X_{i+1}$ & $\left(
i=2,3,5\right) ,$ & $\widehat{\mu }_{2}\left( X_{3},X_{2}\right) =-X_{5},$ &
$\widehat{\mu }_{2}\left( X_{2},X_{4}\right) =X_{6},$ \\
$\widehat{\mu }_{2}\left( X_{7},X_{i}\right) =X_{i}$ & $1\leq i\leq 6,$ & $%
\widehat{\mu }_{2}\left( X_{8},X_{i}\right) =X_{i},$ & $\left(
i=2,3,4\right) ,$ \\
$\widehat{\mu }_{2}\left( X_{8},X_{i}\right) =2X_{i}$ & $\left( i=5,6\right)
.$ &  &
\end{tabular}
\end{equation*}
sur une base $\left\{ X_{1},..,X_{8}\right\} $. On observe que le nilradical de $\widehat{\frak{g}}_{i}$ est isomorphe \`a $\frak{g}_{i}$ pour $i=1,2$. Pour la deuxi\`eme alg\`ebre, nous utilisons une base de vecteurs propres de $\frak{g}_{2}$.

\begin{lemma}
Les alg\`{e}bres $\widehat{\frak{g}}_{1}\otimes\mathbb{C} \simeq \widehat{\frak{g}}_{2}\otimes\mathbb{C} $ sont isomorphes. De plus, l'alg\`{e}bre $\widehat{\frak{g}}_{1}\otimes\mathbb{C}$ est rigide \`{a} cohomologie nulle.
\end{lemma}

En effet, l'alg\`ebre $\widehat{\frak{g}}_{1}\otimes\mathbb{C}$ co\"{\i}ncide avec l'alg\`ebre r\'esoluble rigide $\frak{r}_{8}^{22}$ de la liste \cite{A}. En particulier, il en r\'esulte que  $H^{2}(\frak{g}_{1}\otimes\mathbb{C},\frak{g}_{1}\otimes\mathbb{C})=0$.

\bigskip

\begin{theorem}
Les alg\`ebres de Lie $\widehat{\frak{g}}_{1}$ et $\widehat{\frak{g}}_{2}$ sont rigides compl\`etes sur $\mathbb{R}$. De plus, $\widehat{\frak{g}}_{1} $ n'est pas compl\`{e}tement r\'esoluble.
\end{theorem}

\begin{proof}
D'apr\`{e}s le lemme 3, l'alg\`{e}bre complexifi\'{e}e  $\widehat{\frak{g}}_{2}\otimes\mathbb{C}$ de $\widehat{\frak{g}}_{2} $ est rigide et satisfait $H^{2}(\widehat{\frak{g}}_{2}\otimes\mathbb{C},\widehat{\frak{g}}_{2}\otimes\mathbb{C})=0$. D'apr\`es les propri\'et\'es de la cohomologie, le deuxi\`{e}me
groupe de cohomologie de   $\widehat{\frak{g}}_{1}  $ et $\widehat{\frak{g}}_{1} $ est aussi nul, ce qui montre la rigidit\'{e} en applicant le th\'{e}or\`{e}me de rigidit\'{e} de Nijenhuis et Richardson. Consid\'{e}rons \ l'op\'{e}%
rateur adjoint $ad_{\widehat{\mu }_{1}}\left( X_{8}\right) $. Sur la base $%
\left\{ X_{1},..,X_{8}\right\} $, cet op\'{e}rateur est donn\'{e} par la
matrice
\begin{equation*}
ad_{\widehat{\mu }_{1}}\left( X_{8}\right) =\left(
\begin{array}{cccccccc}
0 & -1 & 0 & 0 & 0 & 0 & 0 & 0 \\
1 & 0 & 0 & 0 & 0 & 0 & 0 & 0 \\
0 & 0 & 0 & 0 & 0 & 0 & 0 & 0 \\
0 & 0 & 0 & 0 & 0 & 1 & 0 & 0 \\
0 & 0 & 0 & 0 & 0 & 0 & 0 & 0 \\
0 & 0 & 0 & -1 & 0 & 0 & 0 & 0 \\
0 & 0 & 0 & 0 & 0 & 0 & 0 & 0 \\
0 & 0 & 0 & 0 & 0 & 0 & 0 & 0
\end{array}
\right) ,
\end{equation*}
et comme le polyn\^{o}me minimal est $m\left( T\right) =T+T^{3}$, $ad_{\widehat{\mu }_{1}}\left( X_{8}\right) $ n'est pas diagonalisable sur $\mathbb{R}$, alors l'alg\`ebre n'est pas compl\'etement r\'esoluble.\newline D'un autre c\^{o}t\'e, l'alg\`ebre $\widehat{\frak{g}}_{2}$ et l'alg\`ebre $r_{8}^{22}$ dans \cite{A} ont les m\^{e}mes constantes de structure, ce qui prouve que cette forme r\'eelle est la normale. La compl\`etude est une cons\'equence de la suite bien connue \cite{Ca}:
\begin{equation*}
\dim Der(\frak{g})=\dim\frak{g}+\dim H^{2}(\frak{g},\frak{g})
\end{equation*}
et la nullit\'e du groupe de cohomologie.
\end{proof}

\medskip
L'alg\`ebre $ \widehat{\frak{g}}_{1}$ admet une d\'ecomposition
$\widehat{\frak{g}}_{1}=\frak{g}_{1}\overrightarrow{\oplus}\frak{t}$,
o\`u $\frak{t}$ est une sous-alg\`ebre ab\'elienne de $Der\frak{g}_{1}$ g\'en\'er\'ee par $X_{7}$ et $X_{8}$. Cependant,
tout vecteur de $\widehat{\frak{g}}_{1} $
qui soit ad-semi-simple fait intervenir une combinaison lin\'eaire
des vecteurs $X_{7}$ et $X_{8}$, et comme  $ad_{\widehat{\mu
}_{1}}\left( X_{8}\right)$ n'est pas diagonalisable, la seule
possibilit\'e c'est consid\'erer des multiples du vecteur $X_{7}$.
De cette fa\c{c}on, au contraire que dans le cas complexe,
l'alg\`ebre $\widehat{\frak{g}}_{1}$ n'admet pas une
d\'ecomposition du type
$\widehat{\frak{g}}_{1}=\frak{g}_{1}\oplus\frak{k}$, o\`{u}
$\frak{k}$ est une sous-alg\`ebre ab\'elienne de $Der(\frak{g}_{1})$
g\'en\'er\'ee par des \'el\'ements $ad$-semi-simples. Dans ce cas
l'alg\`ebre $\frak{t}$ joue le r\^{o}le du tore maximal de
d\'erivations et tout syst\`eme lin\'eaire de racines associ\'e au
vecteur r\'egulier $X_{7}$ est de corang $1$ et  $\dim\frak{t}=2$.
Cette remarque montre que, pour couvrir le cas r\'eel, la
th\'eorie des syst\`emes lin\'eaires obtenue en \cite{A2}, ainsi comme les notions du graphe des poids \cite{C19,C29}, doient
\^{e}tre modifi\'ees. En particulier, on obtient le crit\`ere
suivant:

\begin{proposition}
Soit $\frak{g}$ une alg\`ebre de Lie r\'esoluble rigide complexe.
Alors seulement la forme r\'eelle normale $\frak{g}_{\mathbb{R}}$
est compl\`etement r\'esoluble.
\end{proposition}

\begin{proof}
En effet, si $\frak{g}$ admet des formes r\'eelles, la forme normale a les m\^{e}mes constantes de structure, alors elle est compl\`etement r\'esoluble. Comme en tout cas la partie nilpotente a des constantes de structure sur $\mathbb{R}$, si la forme r\'eelle $\frak{g}^{\prime}$ n'est pas la normale, alors son nilradical poss\`ede une d\'erivation non diagonale sur $\mathbb{R}$ mais sur $\mathbb{C}$ qui appartient au tore ext\'erieur. En cons\'equence, $\frak{g}^{\prime}$ n'est pas compl\`etement r\'esoluble.

\end{proof}

\section{Classification r\'eelle jusqu'\`a dimension 8}

Soit $\frak{g}$ alg\`{e}bre de Lie r\'{e}elle. Alors on a $\dim
H^{2}\left(
\frak{g},\frak{g}\right)  =\dim H^{2}\left(  \frak{g\otimes}\mathbb{C,}%
\frak{g\otimes}\mathbb{C}\right)  $, \'{e}tant donn\'{e} que les
coefficients du syst\`{e}me qui donne les cocycles est \`{a}
coefficients r\'eels. Alors
on peut \'{e}tablir l'\'equivalence suivante:%
\[
H^{2}\left(  \frak{g},\frak{g}\right)  =0\Leftrightarrow
H^{2}\left(
\frak{g\otimes}\mathbb{C,}\frak{g\otimes}\mathbb{C}\right)  =0.
\]
Comme cons\'equence, pour les alg\`{e}bres alg\'{e}briquement
rigides\footnote{On appelle une alg\`ebre alg\'ebriquement rigide si elle satisfait le crit\`ere de Nijenhuis et Richardson (voir \cite{NR}).}, une alg\`{e}bre qui ne soit pas purement complexe est
rigide si et seulement si les formes r\'{e}elles sont rigides. On
peut prouver facilement que pour dimension $n\leq8$, toute
alg\`{e}bre de Lie rigide sur $\mathbb{C}$ est \`{a} cohomologie
nulle. Alors, pour obtenir une classification r\'{e}elle, il
suffit d'ajouter les formes r\'{e}elles correspondantes. Comme en
dimension $n\leq 5$ la classification des alg\`ebres de Lie
nilpotentes sur $\mathbb{R}$  et $\mathbb{C}$
 co\"{\i}ncide, la classification des alg\`ebres de Lie r\'esolubles
 r\'eelles en dimension $n\leq 7$ est la m\^{e}me \cite{A}. L'exemple ant\'erieur
 montre que en dimension 8, la classification r\'eelle diff\`ere de la complexe.

\begin{proposition}
Soit $\frak{g}$ une alg\`{e}bre de Lie r\'{e}soluble alg\'ebriquement rigide de dimension 8
r\'{e}elle dont le tore contient au moins une d\'{e}rivation
non-diagonale. Alors le nilradical $\frak{n}$ est de dimension
$m\leq 6$.
\end{proposition}

\begin{proof}
Soit $\frak{g}$ de rang 1, c'est-\`{a}-dire, $\dim\frak{n}=7$.
Comme tout espace compl\'ementaire de $\frak{n}$ dans $\frak{g}$  est de dimension 1, g\'{e}n\'{e}r\'e par un
vecteur $X$, et l'op\'{e}rateur $ad\left(  X\right)  $ n'est pas
diagonalisable sur $\mathbb{R}$ (mais sur $\mathbb{C)}$, les
valeurs propres de $ad\left( X\right)  $ sont purement imaginaires
ou z\'{e}ro. En effet, la forme de Jordan r\'{e}elle de
$ad\left(  X\right)  $ restreint au nilradical est une matrice du type%
\[
ad\left(  X\right)  \sim\left(
\begin{array}
[c]{cccccccc}%
\lambda_{1} & a_{1} &  &  &  &  &  & \\
-a_{1} & \lambda_{1} &  &  &  &  &  & \\
&  & \ddots &  &  &  &  & \\
&  &  & \lambda_{k} & a_{k} &  &  & \\
&  &  & -a_{k} & \lambda_{k} &  &  & \\
&  &  &  &  & b_{1} &  & \\
&  &  &  &  &  & \ddots & \\
&  &  &  &  &  &  & b_{s}
\end{array}
\right),
\]
o\`u $k\geq 1$ et $\lambda_{i},a_{i},b_{i}\in\mathbb{R}$ . Par la dimension impaire de $\frak{n}$ il existe au moins un valeur propre r\'eel. D'apr\`es la preuve de la proposition 2, les valeurs propres complexes de $ad(X)$ sont un multiple (complexe) d'une suite des nombres enti\`eres. En cons\'equence, ils sont tous r\'eels ou imaginaires pures, d'o\`u $\lambda=0$ est n\'{e}cessairement un valeur propre
du op\'{e}rateur $ad\left(  X\right)  $. Alors l'alg\`{e}bre
rigide complexifi\'{e}e admet z\'{e}ro comme poid. D'apr\`{e}s la
classification de \cite{A2}, une alg\`{e}bre de Lie r\'{e}soluble
de rang un admettant $\lambda=0$ comme poid ne peut pas \^{e}tre rigide,
d'o\`{u} l'affirmation.
\end{proof}

\begin{proposition}
Soit $\frak{g}$ une alg\`ebre de Lie r\'eelle r\'esoluble et alg\`ebriquement rigide de dimension 8. Alors $\frak{g}$ est isomorphe \`a l'une des formes
r\'eelles normales $(\frak{r}_{8,i})_{\mathbb{R}}$ des alg\`ebres
de la liste \cite{A2} o\`u \`a l'une des formes r\'eelles
suivantes:

\begin{itemize}

\item $\frak{r}_{8,22}^{2}=\frak{\widehat{g}}_{1}$:

\begin{equation*}
\begin{tabular}{llll}
$\left[ X_{1},X_{i}\right] =X_{i+1}$ & $\left( 2\leq i\leq
4\right) ,$ & $\left[X_{3},X_{2}\right] =X_{6},$ & $%
\left[ X_{6},X_{2}\right] =X_{5},$ \\
$\left[ X_{7},X_{i}\right] =X_{i}$ & $\left( i=1,2\right) ,
$ & $\left[ X_{7},X_{3}\right] =2X_{3},$ & $\left[ X_{7},X_{5}\right] =4X_{5},$ \\
$\left[ X_{7},X_{i}\right] =3X_{i}$ & $\left( i=4,6\right)
,$ & $\left[ X_{8},X_{1}\right] =X_{2},$ & $\left[ X_{8},X_{2}\right] =-X_{1},$ \\
$\left[ X_{8},X_{4}\right] =-X_{6},$ & $\left[ X_{8},X_{6}\right] =X_{4}.$ &  &
\end{tabular}
\end{equation*}

\item $\frak{r}_{8,29}^{2}$:

\begin{equation*}
\begin{tabular}{llll}
$\left[ X_{1},X_{i}\right] =X_{i+2}$ & $\left( 2\leq i\leq
4\right) ,$ & $\left[X_{2},X_{3}\right] =X_{6},$ & $%
\left[ X_{2},X_{4}\right] =-X_{5},$ \\
$\left[ X_{7},X_{i}\right] =X_{i}$ & $\left( i=1,2\right) ,
$ & $\left[ X_{7},X_{i}\right] =2X_{i},$ & $\left(i=3,4\right),$ \\
$\left[ X_{7},X_{i}\right] =3X_{i}$ & $\left( i=5,6\right)
,$ & $\left[ X_{8},X_{1}\right] =X_{2},$ & $\left[ X_{8},X_{2}\right] =-X_{1},$ \\
$\left[ X_{8},X_{5}\right] =X_{6},$ & $\left[ X_{8},X_{6}\right] =-X_{5}.$ &  &
\end{tabular}
\end{equation*}

\end{itemize}

\end{proposition}

\begin{proof}
Comme les alg\`ebres r\'esolubles rigides complexes de dimension $n\leq 8$ sont rationelles, alors la forme normale des classes d'isomorphisme $\frak{r}_{8,i}$ de la liste \cite{A2} sont rigides comme alg\`ebres r\'eelles. D'apr\`es la remarque ant\'erieure, il suffit d\'eterminer les formes normales des alg\`ebres nilpotentes qui apparaissent comme nilradical d'une alg\`ebre de Lie r\'esoluble complexe rigide en dimension 8 et  poss\`edent des d\'erivations non diagonales sur $\mathbb{R}$. \c{C}a implique que la dimension du nilradical est 6 ou 7, mais par la proposition 5, la dimension ne peut pas \^{e}tre 7. Le probl\`eme est alors r\'eduit \`a voir les formes r\'eelles des alg\`ebres nilpotentes complexes de dimension 6. Parmi ces alg\`{e}bres, seulement $\mathcal{N}_{6,6}\otimes\mathbb{C\simeq
}\mathcal{N}_{6,7}\otimes\mathbb{C}$ et $\mathcal{N}_{6,12}\otimes
\mathbb{C\simeq}\mathcal{N}_{6,14}\otimes\mathbb{C}$ aparaissant comme des
nilradicaux des alg\`{e}bres de Lie rigides complexes $\frak{r}_{8,22}$ et
$r_{8,29}$ de la liste \cite{A}, o\`{u} les alg\`{e}bres $\mathcal{N}%
_{6,6},\mathcal{N}_{6,7}\mathbb{,}\mathcal{N}_{6,12}$ et $\mathcal{N}%
_{6,14}\mathbb{\ }$\ sont d\'{e}fnies par les crochets:

\begin{enumerate}
\item $\mathcal{N}_{6,12}:$%
\[
\left[  X_{1},X_{2}\right]  =X_{4},\;\left[  X_{1},X_{3}\right]
=X_{5},\;\left[  X_{1},X_{4}\right]  =X_{5},\;\left[  X_{2},X_{4}\right]
=X_{6}.
\]

\item $\mathcal{N}_{6,14}:$%
\begin{align*}
\left[  X_{1},X_{2}\right]    & =X_{4},\;\left[  X_{1},X_{3}\right]
=X_{5},\;\left[  X_{1},X_{4}\right]  =X_{6},\;\left[  X_{2},X_{3}\right]
=X_{6},\\
\left[  X_{2},X_{4}\right]    & =-X_{5}.
\end{align*}
\end{enumerate}
Les alg\`{e}bres $\mathcal{N}_{6,6}$ et $\mathcal{N}_{6,7}$ sont celles
employ\'{e}es dans l'exemple de la section 2. L'alg\`{e}bre $\mathcal{N}%
_{6,12}$ admet un tore ext\'{e}rieur $\frak{t}$ de dimension 2, et la somme
semi-directe $\mathcal{N}_{6,12}\overrightarrow{\oplus}\frak{t}$ est isomorphe
\`{a} l'alg\`{e}bre $\frak{r}_{8,29}$ de \cite{A}. L'alg\`{e}bre $\mathcal{N}_{6,14}$
admet seulement une d\'{e}rivation diagonalisable, et une nondiagonalisable.
La somme semi-directe correspond \`{a} $\frak{r}_{8,29}^{2}$.

\end{proof}

\begin{remafin}
On peut se poser la question si pour les alg\`ebres de Lie complexes et r\'eelles la rigidit\'e se preverse par complexification et r\'ealification (comme pour la semi-simplicit\'e). Ici il y a deux probl\`emes: D'abord, on peut construire des alg\`ebres r\'esolubles complexes sans formes r\'eelles (voir \cite{A3}), et la nullit\'e de la cohomologie n'est pas une condition n\'ecessaire. Est-ce que on peut trouver des alg\`ebres r\'eelles rigides $\frak{g}$ sur $\mathbb{R}$ (et dont $\dim H^{2}(\frak{g},\frak{g})\neq 0$ qui ne soient pas rigides sur $\mathbb{C}$? On peut conj\'eturer qu'une alg\`ebre de Lie r\'eelle est rigide (sur le corps r\'eel) si et seulement l'alg\`ebre complexifi\'ee l'est sur le corps complexe. Mais pour les alg\`ebres avec cohomologie non nulle il nous manque d'une m\'ethode effective qui nous permet d'\'etudier la rigidit\'e. Surtout pour rank un, l'existance de formes r\'eelles non normales \'equivaut a la conjecture pos\'ee dans \cite{Cc} concernant l'impossibilit\'e du poid $\lambda=0$ dans le cas complexe.

\end{remafin}

\end{document}